\documentclass[a4paper]{article}
\begin{document} 
\newtheorem{prop}{Proposition}[section]
\newtheorem{Def}{Definition}[section] \newtheorem{theorem}{Theorem}[section]
\newtheorem{lemma}{Lemma}[section] \newtheorem{Cor}{Corollary}[section]

\title{\bf Unconditional well-posedness for the Dirac - Klein - Gordon system 
in two space dimensions}
\author{
{\bf Hartmut Pecher}\\
Fachbereich Mathematik und Naturwissenschaften\\
Bergische Universit\"at Wuppertal\\
Gau{\ss}str.  20\\
42097 Wuppertal\\
Germany\\
e-mail {\tt Hartmut.Pecher@math.uni-wuppertal.de}}
\date{}
\maketitle

\begin{abstract}
The solution of the Dirac - Klein - Gordon system in two space dimensions with 
Dirac data in $H^s$ and wave data in $H^{s+\frac{1}{2}} \times 
H^{s-\frac{1}{2}}$ is uniquely determined in the natural solution space 
$C^0([0,T],H^s) \times C^0([0,T],H^{s+\frac{1}{2}})$, provided $ s > 1/30$ .
This improves the uniqueness part of the global well-posedness result by A. 
Gr\"unrock and the author, where uniqueness was proven in 
(smaller) spaces of Bourgain type. Local well-posedness is also proven for 
Dirac data in $L^2$ and wave data in $H^{\frac{3}{5}+} \times H^{-\frac{2}{5}+} 
$ in the solution space $C^0([0,T],L^2) \times C^0([0,T],H^{\frac{3}{5}+})$ and 
also for more regular data.
\end{abstract}

\renewcommand{\thefootnote}{\fnsymbol{footnote}}
\footnotetext{\hspace{-1.8em}{\it 2000 Mathematics Subject Classification:} 
35Q55, 35L70 }
\normalsize 
\setcounter{section}{0}
\section{Introduction and main results}
The Cauchy problem for the Dirac -- Klein -- Gordon 
equations in two space dimensions reads as follows 
\begin{eqnarray}
\label{0.1}
i(\partial_t + \alpha \cdot \nabla) \psi + M \beta \psi & = & - \phi \beta \psi 
\\
\label{0.2}
(-\partial_t^2 + \Delta) \phi + m\phi & = & - \langle \beta \psi,\psi \rangle
\end{eqnarray}
with (large) initial data
\begin{equation}
\psi(0)  =  \psi_0 \,  , \, \phi(0)  =  \phi_0 \, , \, \partial_t 
\phi(0) = \phi_1 \, .
\label{0.3}
\end{equation}
Here $\psi$ is a two-spinor field, i.e. $\psi : {\bf R}^{1+2} \to {\bf C}^2$, 
and 
$\phi$ is a real-valued function, i.e. $\phi : {\bf R}^{1+2} \to {\bf R}$ , 
$m,M 
\in {\bf R}$ and $\nabla = (\partial_{x_1} , \partial_{x_2}) $ , $ \alpha \cdot 
\nabla = \alpha^1 \partial_{x_1} + \alpha^2 \partial_{x_2}$ . 
$\alpha^1,\alpha^2, \beta$ are hermitian ($ 2 \times 2$)-matrices satisfying 
$\beta^2 = 
(\alpha^1)^2 = (\alpha^2)^2 = I $ , $ \alpha^j \beta + \beta \alpha^j = 0, $  $ 
\alpha^j \alpha^k + \alpha^k \alpha^j = 2 \delta^{jk} I $ . \\
$\langle \cdot,\cdot \rangle $ denotes the ${\bf C}^2$ - scalar product. A 
particular representation is given by \\ $\alpha^1 = {0\;\,1 \choose 1\;\,0}$ , 
$\alpha^2 =  {0\,-i \choose i\;\,0}$ , $\beta = {1\;\,0\choose0 -1}$.\\
We consider Cauchy data in Sobolev spaces: $\psi_0 \in H^s $ , $ \phi_0 \in 
H^r 
$ , $ \phi_1 \in H^{r-1}$ .\\
Local well-posedness was shown by d'Ancona, Foschi and Selberg \cite{AFS1} in 
the case $ s > - \frac{1}{5}$ and $ 
\max(\frac{1}{4}-\frac{s}{2},\frac{1}{4}+\frac{s}{2},s) < r < 
\min(\frac{3}{4}+2s,\frac{3}{4}+\frac{3s}{2},1+s)$. As usually they apply the 
contraction mapping principle to the system of integral equations belonging to 
the problem above. The fixed point is constructed in spaces of Bourgain type 
$X^{s,b} \times X^{r,b}$ which are subsets of the space $C^0([0,T],H^s({\bf 
R}^2)) \times 
C^0([0,T],H^{r}({\bf R}^2))$. Thus especially uniqueness is shown 
also in these spaces of $X^{s,b}$-type.
Thus the question arises whether 
unconditional 
uniqueness holds, namely uniqueness in the natural solution space 
$C^0([0,T],H^s({\bf R}^2)) \times C^0([0,T],H^{r}({\bf R}^2))$ 
without assuming that the solution belongs to some (smaller) $X^{s,b} \times 
X^{r,b}$-space. 

The question of global well-posedness for the system 
(\ref{0.1}),(\ref{0.2}),(\ref{0.3}) was recently answered 
positively for data $\psi_0 \in H^s $ , $ \phi_0 \in 
H^{s+\frac{1}{2}} 
$ , $ \phi_1 \in H^{s-\frac{1}{2}}$ in the case $s \ge 0$ by A. Gr\"unrock and 
the author \cite{GP}. They showed 
existence and uniqueness in Bourgain type spaces $X^{s,b,1}$ based on certain 
Besov spaces with respect to time. These solutions were shown to belong 
automatically to $C^0([0,T],H^s({\bf R}^2)) \times 
C^0([0,T],H^{s+\frac{1}{2}}({\bf R}^2))$. Again the question arises whether 
unconditional 
uniqueness holds, namely uniqueness in the natural solution space 
$C^0([0,T],H^s({\bf R}^2)) \times C^0([0,T],H^{s+\frac{1}{2}}({\bf R}^2))$ 
without assuming that the solution belongs to some (smaller) Bourgain type 
spaces.

The question of unconditional uniqueness was considered among others by Yi Zhou 
for the KdV equation \cite{Z} and nonlinear wave equations \cite{Z1}, by N. 
Masmoudi and K. Nakanishi for the Maxwell-Dirac, the 
Maxwell-Klein-Gordon equations \cite{MN}, the Klein-Gordon-Zakharov system and 
the Zakharov system \cite{MN1}, and by F. Planchon \cite{P} for semilinear wave 
equations.

Our main results read as follows:

\begin{theorem}
\label{Theorem 1.0}
Let $\psi_0 \in H^s({\bf R}^2) $ , $ \phi_0 \in H^r({\bf R}^2) $ , $\phi_1 \in 
H^{r-1}({\bf R}^2) $ , where 
$$\frac{1}{8} >s \ge 0 \quad , \quad \frac{3}{5}-2s < r < 
\min(\frac{3}{4}+\frac{3}{2}s,1-2s) \, .$$ 
Then the Cauchy 
problem (\ref{0.1}),(\ref{0.2}),(\ref{0.3}) is unconditionally locally 
well-posed in
$$ (\psi,\phi,\phi_t) \in C^0([0,T],H^s({\bf R}^2)) \times C^0([0,T],H^r({\bf 
R}^2))  \times C^0([0,T],H^{r-1}({\bf R}^2)) \, . $$
Especially we can choose $ s=0 $ and $ r=\frac{3}{5}+ $.
\end{theorem}
{\bf Remark:} Similar results for $ s \ge \frac{1}{8} $ and a suitable range for 
$r$ can also be given. If $ 1 > s \ge \frac{1}{8} $ the result remains true if
$ \max(\frac{1}{4}+\frac{s}{2},s,\frac{2}{5}-\frac{2}{5}s) < r < 
\min(\frac{3}{4}+\frac{3}{2}s,6s,1) ,$ e.g. if $s = \frac{1}{6}$ and $ 
\frac{1}{3} < r < 1 $.
\\

\begin{theorem}
\label{Theorem 1.1}
Let $\psi_0 \in H^s({\bf R}^2)$ , $ \phi_0 \in H^{s+\frac{1}{2}}({\bf R}^2)$ , 
$ \phi_1 \in H^{s-\frac{1}{2}}({\bf R}^2)$ with $ s > \frac{1}{30}  $. 
Then the Cauchy 
problem (\ref{0.1}),(\ref{0.2}),(\ref{0.3}) is unconditionally globally 
well-posed in 
the space
$$(\psi,\phi,\phi_t) \in C^0({\bf R}^+,H^s({\bf R}^2)) \times 
C^0({\bf R}^+,H^{s+{\frac{1}{2}}}({\bf R}^2)) \times 
C^0({\bf R}^+,H^{s-{\frac{1}{2}}}({\bf R}^2)) \, .$$
This means that existence and uniqueness holds in these spaces.
\end{theorem}
{\bf Remark:} The interesting question of unconditional uniqueness in the case 
of lowest regularity of the data where global existence is known ($s=0$ in 
Theorem \ref{Theorem 1.1} and $s=0 \, , \, r=\frac{1}{2}$ in Theorem 
\ref{Theorem 1.0})(cf. Theorem \ref{Theorem 1.3}) unfortunately remains 
unsolved.

We use the following Bourgain type function spaces. Let $\, \tilde{} \,$ denote 
the Fourier transform with respect to space and time. $X^{s,b}_{\pm}$ is the 
completion of ${\cal S}({\bf R} \times {\bf R}^2)$ with respect to
$$ \|f\|_{X_{\pm}^{s,b}} = \|U_{\pm}(-t)f\|_{H^b_t H^s_x} = \| \langle \xi 
\rangle^s \langle \tau \pm |\xi| \rangle^b \tilde{f}(\xi,\tau)\|_{L^2} \, , $$
where $U_{\pm}(t) := e^{\mp it|D|}$ and 
$$ \|g\|_{H^b_t H^s_x} = \| \langle \xi \rangle^s \langle \tau \rangle^b 
\tilde{g}(\xi,\tau)\|_{L^2_{\xi \tau}} \, . $$
Finally we define
$$ \|f\|_{X^{s,b}_{\pm}[0,T]} := \inf_{g_{|[0,T]}=f}  \|g\|_{X^{s,b}_{\pm}} \, 
. $$

\section{Preparations}
As is well-known it is convenient to replace the system 
(\ref{0.1}),(\ref{0.2}),(\ref{0.3}) by 
considering the projections onto the one-dimensional eigenspaces of the 
operator 
$-i \alpha \cdot \nabla$ belonging to the eigenvalues $ \pm |\xi|$. These 
projections are given by $\Pi_{\pm}(D)$, where  $ D = 
\frac{\nabla}{i} $ and $\Pi_{\pm}(\xi) = \frac{1}{2}(I 
\pm \frac{\xi}{|\xi|} \cdot \alpha) $. Then $ 
-i\alpha \cdot \nabla = |D| \Pi_+(D) - |D| \Pi_-(D) $ and $ \Pi_{\pm}(\xi) 
\beta
= \beta \Pi_{\mp}(\xi) $. Defining $ \psi_{\pm} := \Pi_{\pm}(D) \psi$ and 
splitting the function $\phi$ into the sum $\phi = \frac{1}{2}(\phi_+ + 
\phi_-)$, where $\phi_{\pm} := \phi \pm iA^{-1/2} \partial_t \phi $ , $ A:= 
-\Delta+1$ , the Dirac - Klein - Gordon system can be rewritten as
\begin{eqnarray}
\label{*}
(-i \partial_t \pm |D|)\psi_{\pm} & = & -M\beta \psi_{\mp} + \Pi_{\pm}(\phi 
\beta (\psi_+ + \psi_-)) \\
\nonumber
(i\partial_t \mp A^{1/2})\phi_{\pm} & = &\mp A^{-1/2} \langle \beta (\psi_+ + 
\psi_-), \psi_+ + \psi_- \rangle \mp A^{-1/2} (m+1)(\phi_+ + \phi_-) . \\
\label {**}
\end{eqnarray}
The initial conditions are transformed into
\begin{equation}
\label{***}
\psi_{\pm}(0) = \Pi_{\pm}(D)\psi_0 \, , \,  \phi_{\pm}(0) = \phi_0 \pm i 
A^{-1/2} \phi_1
\end{equation}

We now state again the above mentioned well-posedness results on which our 
results rely.
\begin{theorem}
\label{Theorem 1.2} (\cite{AFS1})
Let $\psi_0 \in H^s $ , $ \phi_0 \in H^r$ , $ \phi_1 \in H^{r-1} $ , where
$$ s > - \frac{1}{5} \,\, , \, \, 
\max(\frac{1}{4}-\frac{s}{2},\frac{1}{4}+\frac{s}{2},s) < r < 
\min(\frac{3}{4}+2s,\frac{3}{4}+\frac{3}{2}s,1+s) \, . $$
Then the Cauchy problem (\ref{*}),(\ref{**}),(\ref{***}) is locally well-posed 
for
$$ (\psi_{\pm},\phi_{\pm}) \in X^{s,\frac{1}{2}+}_{\pm}[0,T] \times 
X^{r,\frac{1}{2}+}_{\pm}[0,T]\, , $$
i.e.
\begin{eqnarray*}
(\psi,\phi,\partial_t \phi) & \in & (X^{s,\frac{1}{2}+}_+[0,T] + 
X^{s,\frac{1}{2}+}_-[0,T]) 
\times (X^{r,\frac{1}{2}+}_+[0,T] + X^{r,\frac{1}{2}+}_-[0,T]) \\
& & \times 
(X^{r-1,\frac{1}{2}+}_+[0,T] + X^{r-1,\frac{1}{2}+}_-[0,T]) \, . 
\end{eqnarray*}
This solution belongs to
$$ C^0([0,T],H^s) \times  C^0([0,T],H^r) \times C^0([0,T],H^{r-1}) \, . $$
\end{theorem}
{\bf Remark:} The question of uniqueness in the latter (larger) spaces remained 
open.
\begin{theorem}
\label{Theorem 1.3} (\cite{GP})
Let $ s \ge 0$ and $ \psi_0 \in H^s $ , $ \phi_0 \in H^{s+\frac{1}{2}} $ , $ 
\phi_1 \in H^{s-\frac{1}{2}} $ . Then the Cauchy problem 
(\ref{*}),(\ref{**}),(\ref{***}) is globally well-posed for 
$$ (\psi_{\pm},\phi_{\pm}) \in X^{s,\frac{1}{3},1}_{\pm} \times 
X^{s+\frac{1}{2},\frac{1}{3},1}_{\pm} \, . $$
This solution belongs to
$$ (\psi,\phi,\partial_t \phi) \in C^0({\bf R}^+,H^s) \times C^0({\bf 
R}^+,H^{s+\frac{1}{2}}) \times C^0({\bf R}^+,H^{s-\frac{1}{2}}) \, . $$
Here the spaces $X^{s,\frac{1}{3},1}$ are certain Bourgain type spaces based on 
Besov spaces (with respect to time). For a precise definition we refer to 
\cite{GP}.
\end{theorem}
{\bf Remark:} Again the question of uniqueness in the latter (larger) spaces 
remained open.

We recall the following facts about the solution of the inhomogeneous linear 
problem
$$ \partial_t v - i\phi(D) v = F \quad , \quad v(0) = v_0 \, , $$
namely 
$$ v(t) = U(t)v_0 + \int_0^t U(t-s) F(s) ds \, , $$
where
$$ U(t) = e^{it\phi(D)} v_0 \, . $$
\begin{prop}
\label{Prop.} (\cite{GTV} or \cite{G})
Let $ b'+1 \ge b \ge 0 \ge b' > -1/2 $. Then the following estimate holds for 
$T \le 1$:
$$ \|v\|_{X^{s,b}[0,T]} \le c(T^{\frac{1}{2}-b} \|v_0\|_{H^s} + 
T^{1+b'-b}\|F\|_{X^{s,b'}[0,T]}) \, . $$
Here $X^{s,b}$ denotes the completion of ${\cal S}({\bf R} \times {\bf R}^2)$ 
with respect to the norm $\|f\|_{X^{s,b}} = \|U(-t)f\|_{H^b_t H^s_x}$ and 
$X^{s,b}[0,T]$ the restrictions of these functions to $[0,T]$.
\end{prop}

\section{Proofs of the theorems}
The key result reads as follows:
\begin{theorem}
\label{Theorem 2.1}
Let $\psi_0 \in H^s({\bf R}^2)$ , $\phi_0 \in H^r({\bf R}^2)$ , 
$\phi_1 \in H^{r-1}({\bf R}^2) $ , $ T > 0 $ . Assume $ \frac{1}{8} > s \ge 0$ 
and $ \frac{3}{5}-2s < r < 1-2s $.  Then the Cauchy problem 
(\ref{0.1}),(\ref{0.2}),(\ref{0.3}) has at most one solution
$$ (\psi,\phi,\partial_t \phi) \in C^0([0,T],H^s({\bf R}^2)) \times 
C^0([0,T],H^r({\bf R}^2))  \times C^0([0,T],H^{r-1}({\bf 
R}^2)) \, . $$
This solution satisfies $\psi_{\pm} \in 
X_{\pm}^{-\frac{1}{2}+ \frac{r}{2}+s+,\frac{1}{2}+}[0,T]$ , $ \phi_{\pm} \in 
X_{\pm}^{-\frac{1}{4}+r+2s+,\frac{1}{2}+}[0,T]$ .
\end{theorem}
{\bf Proof:} We show that any solution
$$ (\psi,\phi,\partial_t \phi) \in C^0([0,T],H^s({\bf R}^2)) \times 
C^0([0,T],H^r({\bf R}^2))  \times C^0([0,T],H^{r-1}({\bf 
R}^2)) \,  $$
fulfills 
$\psi_{\pm} \in 
X_{\pm}^{-\frac{1}{2}+ \frac{r}{2}+s+,\frac{1}{2}+}[0,T]$ , $ \phi_{\pm} \in 
X_{\pm}^{-\frac{1}{4}+r+2s+,\frac{1}{2}+}[0,T]$ .
In this space uniqueness holds by the 
result of d'Ancona, Foschi and Selberg (Theorem \ref{Theorem 1.2}), who had to 
use the 
full null structure of the system.

Let $\psi_{\pm} \in C^0([0,T],H^s)$ , $ \phi_{\pm} \in 
C^0([0,T],H^r)$ be a solution of (\ref{*}),(\ref{**}),(\ref{***}) 
in the interval [0,T] for some $ T \le 1 $.\\
{\bf a.} We estimate
\begin{eqnarray*}
 \|\phi \beta \psi_{\pm} \|_{L^2((0,T),H^{-1+r+s})}  
& \le & c \|\phi \beta \psi_{\pm}\|_{L^2((0,T),L^{\tilde{r}})} \\
& \le & c T^{\frac{1}{2}} \|\phi\|_{L^{\infty}((0,T),L^{\tilde{p}})} 
\|\psi_{\pm}\|_{L^{\infty}((0,T),L^{\tilde{q}})}  \\
& \le & c T^{\frac{1}{2}} \|\phi\|_{L^{\infty}((0,T),H^{s+\frac{1}{2}})} 
\|\psi_{\pm}\|_{L^{\infty}((0,T),H^s)} < \infty \, ,
\end{eqnarray*}
where $ \frac{1}{\tilde{r}} = 1-\frac{r}{2}-\frac{s}{2}$ , $ \frac{1}{\tilde{p}} 
= 
\frac{1}{2}-\frac{r}{2} $ , $ \frac{1}{\tilde{q}} = \frac{1}{2} - \frac{s}{2} $ 
. \\
We also have $ \psi_{\pm} \in L^2((0,T),H^{-1+r+s}) ,$ because $r<1$, so that 
from 
(\ref{*}) we get
$ \psi_{\pm} \in X^{-1+r+s,1}_{\pm}[0,T] ,$ because
$$ \|\psi_{\pm}\|^2_{X_{\pm}^{-1+r+s,1}[0,T]} \sim \int_0^T 
\|\psi_{\pm}(t)\|^2_{H^{-1+r+s}} dt + \int_0^T \|(-i\partial_t \pm 
|D|)\psi_{\pm}(t)\|^2_{H^{-1+r+s}} ds \, . $$
Interpolation with $\psi_{\pm} \in X_{\pm}^{s,0}[0,T]$ gives $ \psi_{\pm} \in 
X_{\pm}^{s_1,\frac{1}{2}+}[0,T]$ , where $ s_1 = -\frac{1}{2}+\frac{r}{2}+s+  . 
$ Remark that $s_1 <0$ under our assumptions.\\
{\bf b.} In order to show from (\ref{**}) that $\phi_{\pm} \in 
X_{\pm}^{r_1,\frac{1}{2}+}[0,T] $ we have to give the following 
estimates according to Prop. \ref{Prop.}:
\begin{enumerate}
\item $$ \| \langle \beta \Pi_{\pm 1}(D)\psi, \Pi_{\pm 2} \psi ' \rangle 
\|_{X_{\pm 3}^{r_1 -1,-\frac{1}{2}+}[0,T]} \le c \|\psi\|_{X_{\pm 
1}^{s_1,\frac{1}{2}+}[0,T]} \|\psi'\|_{X_{\pm 
2}^{s_1,\frac{1}{2}+}[0,T]} $$
Here $\pm_1,\pm_2,\pm_3$ denote independent signs. This estimate is proven in 
\cite{AFS1}, Thm. 2 and requires the following conditions: $ s_1 > -\frac{1}{4} 
\, \Leftrightarrow \, r+2s > \frac{1}{2} $ and $r_1 < \frac{3}{4} +2s_1 = 
-\frac{1}{4}+r+2s+ $. Thus we can choose $r_1 = -\frac{1}{4}+r+2s+ $ .
\item 
$$\|A^{-\frac{1}{2}} \phi_{\pm}\|_{X_{\pm}^{r_1,-\frac{1}{2}+}[0,T]} 
\le 
\|\phi_{\pm}\|_{L^2((0,T),H^{r_1-1})} \le T^{\frac{1}{2}} 
\|\phi_{\pm}\|_{L^{\infty}((0,T),H^{r_1-1})} < \infty $$
\item $ \phi_{\pm}(0) \in H^r \subset H^{r_1} $ , if $ s < \frac{1}{8} 
$.
\end{enumerate}
Choosing $\psi = \psi_{\pm 1}$ and $\psi ' = \psi_{\pm 2}$ in 1. and using 2. 
and a. we get $ \phi_{\pm} \in X_{\pm}^{r_1,\frac{1}{2}+}[0,T] .$ \\
{\bf c.} We have shown that any solution $\psi_{\pm} \in C^0([0,T],H^s)$ , $ 
\phi_{\pm} \in C^0([0,T],H^r)$  fulfills $ \psi_{\pm} \in X_{\pm}^{s_1}[0,T] $ , 
$ \phi_{\pm} \in X_{\pm}^{r_1}[0,T] $ . Now we use the uniqueness part of 
Theorem \ref{Theorem 1.1}. It requires the following conditions:
$$ \max(\frac{1}{4}-\frac{s_1}{2},\frac{1}{4}+\frac{s_1}{2},s_1) < r_1 
<\min(\frac{3}{4}+2s_1,\frac{3}{4}+\frac{3}{2}s_1,1+s_1) $$
and $ s_1 > -\frac{1}{5} $. An elementary calculation shows that this is 
equivalent to
$$ \frac{3}{5} -2s < r < 1-2s \, . $$
This gives the claimed result.\\[1em]
{\bf Proof of Theorem \ref{Theorem 1.0}} We combine Theorem \ref{Theorem 2.1} 
with the existence part of the local well-posedness result of d'Ancona, Foschi 
and Selberg (Theorem \ref{Theorem 1.2}). One easily checks that the conditons on 
$s$ and $r$ reduce to the assumed ranges for these parameters.\\[1em]
{\bf Proof of Theorem \ref{Theorem 1.1}:} We use Theorem \ref{Theorem 2.1} with 
$ s<\frac{1}{8} \, , \, r = s+\frac{1}{2}$. This requires $\frac{3}{5}-2s < 
s+\frac{1}{2} \Leftrightarrow s > \frac{1}{30}$ . Combining this
with the existence part of the global well-posedness of A. Gr\"unrock and the 
author (Theorem \ref{Theorem 1.3}) we get the claimed result.

\end{document}